\documentclass{elsart}
\usepackage{natbib} 
\usepackage{amsmath,amsfonts,amssymb,times}

\begin{document}
\begin{frontmatter}

\title{Localized large sums of random variables}
\author{Kevin Ford\thanksref{NSF}\corauthref{cor}}
\address{Department of Mathematics, University of Illinois at
   Urbana-Champaign,\goodbreak 1409 West Green St., Urbana, IL, 61801, USA}
\ead{ford@math.uiuc.edu}

\author{G\'erald Tenenbaum}
\address{Institut \'Elie Cartan, Universit\'e Henri--Poincar\'e Nancy 
1,\goodbreak
B.P. 239, 54506 Vand\oe uvre-l\`es-Nancy Cedex, France}
\ead{gerald.tenenbaum@iecn.u-nancy.fr}

\thanks[NSF]{Research supported by National Science
   Foundation Grants DMS-0301083 and DMS-0555367}
\corauth[cor]{Corresponding author.}

\dedicated{Dedicated to the memory of Walter Philipp}
\date{\today}


\begin{abstract}We study large partial sums, localized with respect to the
      sums of
      variances, of a sequence of centered random variables.  An
      application is given to the
      distribution of prime factors of typical integers.
\end{abstract}


\end{frontmatter}

\theoremstyle{plain}


\newcommand{\ZZ}{{\mathbb Z}}
\newcommand{\Q}{{\mathbb Q}}
\newcommand{\R}{{\mathbb R}}
\newcommand{\C}{{\mathbb C}}
\newcommand{\N}{{\mathbb N}}
\newcommand{\PP}{{\mathbb P}}
\newcommand{\EE}{{\mathbb E}}

\newcommand{\be}{\begin{equation}}
\newcommand{\ee}{\end{equation}}
\newcommand{\benn}{\begin{equation*}}
\newcommand{\eenn}{\end{equation*}}

\newcommand{\lm}{\ensuremath{\lambda}}
\newcommand{\lam}{{\ensuremath{\lambda}}}
\newcommand{\Chi}{{\ensuremath{\chi}}}
\renewcommand{\a}{{\ensuremath{\alpha}}}
\renewcommand{\b}{{\ensuremath{\beta}}}
\newcommand{\del}{{\ensuremath{\delta}}}
\newcommand{\eps}{{\ensuremath{\varepsilon}}}
\renewcommand{\(}{\left(}
\renewcommand{\)}{\right)}
\newcommand{\ds}{\displaystyle}
\newcommand{\pfrac}[2]{\left(\frac{#1}{#2}\right)}
\font\fourteengp=cmmi10 at 14pt
    \font\thirteengp=cmmib10 at 13pt
    \font\twelvegp=cmmib10 at 12pt
    \font\elevengp=cmmib10 at 11pt
    \font\tengp=cmmib10
    \font\ninegp=cmmib10 at 9pt
    \font\eightgp=cmmib8
    \font\sevengp=cmmib7
    \font\sixgp=cmmib6
    \font\fivegp=cmmib5

   \def\gp{\textfont1=\tengp\scriptfont1=\sevengp\scriptscriptfont1=\fivegp}

\newcommand{\bo}{\mbox{$\gp\omega$}}
\newcommand{\fl}[1]{\lfloor #1 \rfloor}

\renewcommand{\le}{\leqslant}
\renewcommand{\leq}{\leqslant}
\renewcommand{\ge}{\geqslant}
\renewcommand{\geq}{\geqslant}
\newcommand{\tw}{\widetilde{\varrho}}

\font\tenrsfs=rsfs7 at 10pt
\font\sevenrsfs=rsfs7
\font\fiversfs=rsfs5
\newfam\rsfsfam\textfont\rsfsfam=\tenrsfs\scriptfont\rsfsfam=\sevenrsfs\scriptscriptfont\rsfsfam=\fiversfs
\def\curly{\fam\rsfsfam\tenrsfs}

\def\dater{\vglue-10mm\rightline{(\the\day/\the\month/\the\year)}}
\def\dateheure{\the\day/\the\month/\the\year,\ \heure}
\def\chapter{chap.\thinspace}
\def\pages{pp.\thinspace}
\def\e{{\rm e}}
\def\AA{{\curly A}}
\def\LL{{\curly L}}

\section{Introduction}

Consider random
variables $X_1, X_2, \ldots$  with $\EE X_j=0$ and $\EE X_j^2 = \sigma_j^2$.
Let
$$S_n=X_1+\cdots+X_n, \quad s_n^2 =  \sigma_1^2+\cdots+\sigma_n^2,$$
and assume that (a) $s_n\to \infty$ as $n\to \infty$.\par
Given a positive function $f_N\geqslant 1+1/N$, we are interested in
the behavior~of
$$
I = \liminf_{N\to \infty} \max_{N < s_n^2\le Nf_N} |S_n|/s_n.
$$
If we replace $\liminf$ by $\limsup$, it
immediately follows from the law of the iterated logarithm that
$I=\infty$ almost surely when $f_N$ is bounded.
Our results answer a question
originally raised, in oral form, by A. S\'ark\"ozy
and for which a partial answer had previously
been given by the second
author, see Chap.\thinspace3 of \citet{Oon}.

\section{Independent random variables}

Assume that the $X_j$ are independent.  Then $\EE S_n^2 = s_n^2$.
In addition to condition (a), we will work with two other mild
assumptions, (b) $s_{j+1}/s_j \ll1$ when $s_j>0$
and~(c) for every $\lam>0$, there is a constant $c_\lam>0$
such  that if $n$ is large enough and $s_m^2>2s_n^2$, then
$$
\PP \( | S_m-S_{n}| \ge \lam s_{m} \) \ge c_\lam.
$$
Condition (b) says that no term in $S_n$ dominates the others.  Condition (c)
follows if the Central Limit Theorem (CLT) holds for the sequence of
$S_n$, since
CLT for $S_n$ implies CLT for $S_{m}-S_{n}$ as $(m-n)\to\infty$.  For example,
(c) holds for i.i.d. random variables, under the Lindeberg condition
$$
\forall \eps>0, \quad \lim_{n\to\infty} \sum_{1\leqslant j\leqslant n} \EE \(
X_j^2/s_n^2 : |X_j| > \eps s_n \) = 0
$$
and the stronger Lyapunov condition
$$
\exists \delta>0 : \sum_{1\leqslant j\leqslant n} \EE |X_j|^{2+\del}
= o(s_n^{2+\del}).
$$
Condition (c) is weaker, however, than CLT.

\begin{thm}\label{thm1}
(i) Suppose (a), (b), and $f_N = (\log N)^M$
for some constant $M>0$.  Then $I<\infty$ almost surely.\par
(ii) Suppose (a), (b), (c) and $f_N=(\log N)^{\xi(N)}$ with $\xi(N)$ tending
monotonically to~$\infty$.  Then $I=\infty$ almost surely.
\end{thm}

\noindent{\bf Remark.}
In the first statement of the theorem we show
in fact that almost surely $I\le 15 \sqrt{M+1}(\max_{s_j>0} s_{j+1}/s_j)^2$.

\begin{lem}[Kolmogorov's inequality, 1929]\label{kol}
We have $$\ds \PP(\textstyle\max_{1\leqslant j\leqslant k}|S_j| \ge
\lambda s_k) \le
1/\lambda^2\qquad (k\geqslant 1).$$
\end{lem}

\emph{Proof of Theorem \ref{thm1}}.
By (a) and (b), there is a constant $D$ so that $s_{j+1}/s_j \le D$ for all
large $j$.  Define
$$h(n):=\max\{k:s_k^2\leqslant n\}\qquad (n\in\N^*), $$
so that the conditions $N< s_n^2\leqslant Nf_N$ and $h(N)<n\leqslant
h(Nf_N)$ are equivalent.
\goodbreak\par\goodbreak
We first consider the case when $f_N:=(\log N)^M$.
Let
$$
N_j := j^{(M+3)j}, \quad t(j) := \fl{(M+1)(\log j)/\log 2}, \quad H_j
:= 2^{t(j)},
$$
and
$$
U_j := h(N_j), \quad U_{j,t} := h(2^t N_j) \;\; (0\le t\le t(j)), \quad
V_j: = h(H_j N_j) = U_{j,t(j)}.
$$
It is possible that $U_{j,t+1}=U_{j,t}$ for some $t$.
Note that for large $j$, $H_jN_j \ge N_j f_{N_j}$.
\par\goodbreak
Let $k$ be a constant depending only on $M$ and $D$.
For $j\ge 1$ define the events
\begin{align*}
A_j &:= \{ |S_{V_j}| \le s_{U_{j+1}} \}, \\
B_j &:= \bigcap_{0\le t\le t(j)-1} B_{j,t} \; \text{ where } \;
       B_{j,t}:=
       \left\{ \max_{U_{j+1,t} \le n\le U_{j+1,t+1}} |S_{U_{j+1,t+1}}-S_n|
       \le k s_{U_{j+1,t}} \right\}, \\
C_j &:= \{ |S_{U_{j+1}}-S_{V_j}| \le 2 s_{U_{j+1}} \}.
\end{align*}
By (b) and the definition of $h(N)$, we have
\be\label{sbounds}
D^{-1} \sqrt{2^t N_j} \le s_{U_{j,t}} \le \sqrt{2^t N_j}
\ee
for all $j,t$.
It follows from Lemma \ref{kol} that
$$
\PP(\overline{A_j}) \le D^2 \frac{H_j N_j}{N_{j+1}} \le \frac{D^2}{j^2}.
$$
Thus, $\sum_{j\ge 1} \PP(\overline{A_j}) < \infty$ and hence
almost surely there is a $j_0$ so that $A_j$ occurs for $j\ge j_0$.
Applying Lemma \ref{kol} again yields
$$
\PP(\overline{B_{j,t}}) \le
\frac{s^2_{U_{j+1,t+1}}-s^2_{U_{j+1,t}}}{k^2 s^2_{U_{j+1,t}}} \le
\frac{D^2 2^{t+1} N_{j+1}}{k^2 2^t N_{j+1}} = \frac{2D^2}{k^2}.
$$
If $k=3D\sqrt{M+1}$, then
$$
\PP(B_j) \ge \( 1 - \frac{2D^2}{k^2} \)^{t(j)} \ge \frac{1}{j^{1/2}}
$$
for large $j$.  Also by Lemma \ref{kol}, $\PP(C_j) \ge \frac34$,
and since $B_j$ and $C_j$ are independent,
$$
\sum_{j\ge 1} \PP(B_j C_j) = \infty.
$$
Since the events $B_jC_j$ are independent, the Borel--Cantelli lemma
implies that almost surely the events $B_j C_j$
occur infinitely often.  Thus,  the event $A_jB_jC_j$ occurs for an
infinite sequence of integers $j$.
Take such a index $j$, let $n\in [U_{j+1}, V_{j+1}]$ and
$U_{j+1,g-1} < n \le U_{j+1,g}$, where $1\le g\le t(j+1)$.  We have
by several applications of \eqref{sbounds}
\begin{align*}
|S_n| &\le |S_{V_j}| + |S_{U_{j+1}}-S_{V_j}| + \sum_{0\le t\le g-2}
|S_{U_{j+1,t}}-S_{U_{j+1,t+1}}| + |S_n-S_{U_{j+1,g-1}}| \\
&\le 3s_{U_{j+1}}+k\sum_{0\le t\le g-1} s_{U_{j+1,t}} \\
&\le\big\{3+k(1+2^{1/2}+\cdots+2^{(g-1)/2})\big\}\sqrt{N_{j+1}} \\
&\le 5k \sqrt{2^{g-1} N_{j+1}} \\
&\le 5kD s_n = 15 D^2 (M+1)^{1/2} s_n.
\end{align*}
This completes the proof of part (i) of the theorem, since
$$V_{j+1} \ge h(\tfrac12 j^{M+1} N_j) \ge h(N_j\log^M N_j)$$ for large $j$.

%
%

Now suppose $f_N=(\log N)^{\xi(N)}$ with $\xi(N)$ tending monotonically to
$\infty$. \par
Let $\lam>0$ be arbitrary and define $K:=2D^2$.
Let $N^*_1$ be so large that $f_{N_1^*} \ge K$.
For $j\ge 1$ let $N^*_{j+1}=N^*_j K^{u(j)}$, where
$u(j):=\fl{\log f_{N^*_j}/\log K}$.  Put
$$
U_j^*:=h(N^*_j), \quad U^*_{j,t} := h(K^t N^*_j) \; (0\le t\le u(j)).
$$
Let $J_j:=[U^*_j,U^*_{j+1}]$ and
$$
Y_j := \max_{n\in J_j} |S_n|/s_n.
$$
We have
$$
u(j)\ge 1 \Rightarrow N_{j+1}^*\ge K N_j^* \Rightarrow u(j)/\log j \to \infty.
$$
Therefore, by (c), if $j$ is sufficiently large then
\begin{align*}
\PP(Y_j \le \lam/2 ) &\le \prod_{1\le t\le u(j)} \PP \( |S_{U^*_{j,t}} -
       S_{U^*_{j,t-1}}| \le \tfrac 12
\lambda(s_{U^*_{j,t}} + s_{U^*_{j,t-1}}) \) \\
&\le  \prod_{1\le t\le u(j)} \PP \( |S_{U^*_{j,t}} -  S_{U^*_{j,t-1}}|
\le \lam \sqrt{K^t N_j^*} \) \\
&\le (1-c_\lam)^{u(j)} \le \frac{1}{j^2}.
\end{align*}
Thus
$$
\sum_{j\ge 1} \PP(Y_j\le \lambda/2) < \infty.
$$
Almost surely, $Y_k \le \lam/2$ for only finitely many $k$.

Theorem \ref{thm1} has an analog for Brownian motion, which follows
from Theorem 1 and the invariance principle.

\begin{thm}\label{Brownian}  Let $W(t)$ be Brownian motion on $[0,\infty)$.
If $f_N=(\log N)^M$ with fixed $M>0$, then almost surely
$$
I = \liminf_{N\to \infty} \max_{N< t\le Nf_N}
\frac{|W(t)|}{\sqrt{t}} < \infty.
$$
If $f_N=(\log N)^{\xi(N)}$ with $\xi(N)\to\infty$, then $I=\infty$
almost surely.
\end{thm}

Theorem \ref{Brownian} can be proved directly and more swiftly using the
methods used to establish Theorem \ref{thm1}.  By invariance
principles \citep[e.g.][]{Ph},
   one may deduce from Theorem \ref{Brownian} a version of Theorem
\ref{thm1} where stronger hypotheses on the~$X_j$ are assumed.  As it stands,
now, however, Theorem \ref{thm1} does not follow from Theorem \ref{Brownian}.

%
\section{Dependent random variables}
%

The conclusions of Theorem \ref{thm1} can also be shown to hold for certain
sequences of weakly dependent random variables by making use of almost sure
invariance principles.  We assume that (d) there exists a sequence of
i.i.d. normal random variables $Y_j$ with $\EE Y_j^2=\sigma_j^2$, defined on
the same probability space as the sequence of $X_j$, and such that if
$Z_n=Y_1+\cdots + Y_n$, then
$$
\left| S_n - Z_n \right| = O(s_n) \qquad \mbox{a.s.}
$$
Of course the variables $Y_j$ are dependent on the $X_j$, but not on
each other.  Property (d) has been proved for martingale difference
sequences, sequences satisfying certain mixing conditions, and lacunary
sequences $X_j = \{ n_j \omega \}$ with $\inf n_{j+1}/n_j>1$, $\omega$
uniformly distributed in $[0,1]$ and $\{x\}$ is the fractional part of $x$.
See e.g. \citet{Ph} for a survey of such results.

\begin{thm}\label{thm2}
(i) Suppose (a), (b), and (d).  If $f_N := (\log N)^M$
for some constant $M>0$, then $I<\infty$ almost surely.\par
(ii) Let $\xi(N)$ tend
monotonically to $\infty$ and set $f_N:=(\log N)^{\xi(N)}$. Then
$I=\infty$ almost surely.
\end{thm}

By (d),
$$
I = O(1) + \liminf_{N\to \infty} \; \max_{N< s_n^2 \le Nf_N}
|Z_n|/s_n,
$$
and we apply Theorem \ref{thm1} to the sequence of $Y_j$.  The variable $Z_n$
is normal with variance $s_n^2$, hence (c) holds.

%
\section{Prime factors of typical integers}
%

Consider a sequence of independent random variables $Y_p$, indexed by prime
numbers $p$, such that $\PP(Y_p=1)=1/p$ and $\PP(Y_p=0)=1-1/p$.
We can think of $Y_p$ as modelling whether or not a ``random'' integer
is divisible by $p$.  As $\EE Y_p = 1/p$, we form the centered r.v.'s
$X_p=Y_p-1/p$ (we may also define $X_j$ for non-prime $j$ to be
zero with probability 1).  Let
$$
T_n = \sum_{p\le n} Y_p, \qquad S_n = \sum_{p\le n} X_p.
$$
We have $\EE X_p^2=(1-1/p)/p$, hence by Mertens' estimate
$$
s_n^2 = \sum_{p\le n} \frac{1}{p}-\frac{1}{p^2} = \log_2 n + O(1).
$$
Here and in the sequel, $\log_k$ denotes,  for integer $k\ge 2$, the
$k$-fold iterated logarithm.
Since $\EE |X_p|^3 \le 1/p$, the Lyapunov condition holds with
$\delta=1$.  Then (a), (b) and (c) hold, and therefore the conclusion
of Theorem 1 holds.  Here take $D=\max_{n\ge 2} s_{n+1}/s_n$ since $s_1=0$.

Let $\omega(m,t)$ denote the number of distinct prime factors of $m$
which are $\le t$.
The sequence $\{ T_n  : n\ge 1 \}$ mimics well the behavior of
the function $\omega(m,n)$ for a ``random'' $m$, at least when $n$ is not
too close to $m$.  This is known as the Kubilius
model. It can be made very precise, see \citep[][Ch. 3,
especially \pages119--122]{El} and \citet{GT99} for the sharpest
estimate known to date.  Suppose $r$ is an integer with $2\le r\le x$ and
$r=x^{1/u}$, $\bo_r(m)= (\omega(m,1),\ldots,\omega(m,r))$
and suppose $Q$ is any subset of $\ZZ^r$.
Then, given arbitrary $c<1$, and uniformly in $x,r$ and $Q$, we have
\be\label{Kubilius}
\frac{1}{x} | \{ m\le x : \bo_r(m)\in Q \} | = \PP \bigl(
(T_1,\ldots,T_r) \in Q \bigr) + O\( x^{-c} +
\e^{- u\log u} \).
\ee
An analog of Theorem \ref{thm1}, established by parallel estimates,
provides via \eqref{Kubilius} information about
localized large values~of
$$
\varrho(m,t) := |\omega(m,t)-\log_2 t|/\sqrt{\log_2 t}.
$$

\begin{thm}\label{omega}
(i) Let $M>0$ be fixed, $f_N:=(\log N)^M$ and put $K:=30D^2 \sqrt{M+1}$.  If
$g=g(m)\to\infty$ monotonically as $m\to\infty$ in such a way that
$g^2 f_{g^2} \le \log_2 m$
for large $m$, then for a set of integers $m$ of natural
density $1$,\footnote{A subset $\AA$ of $\N^*$ is said to have
natural density 1 if
$|\AA\cap[1,x]|=x+o(x)$ as $x\to\infty$.} we have
$$
\min_{g(m)\leqslant N\leqslant g(m)^2} \; \max_{N<\log_2 t \le Nf_N}
\varrho(m,t) \le K.
$$
(ii) Let $\xi(N)\to\infty$ in such a way that $f_N:=(\log N)^{\xi(N)}\le N$.
Suppose that $g(m)\to\infty$ monotonically as $m\to\infty$, that
$g(m)\le (\log_2 m)^{1/10}$, and let
$$
I_m := \min_{\substack{g(m) \le N \\ Nf_N \le \log_2 m}} \;\; \max_{N
     \le \log_2 t\le Nf_N} \varrho(m,t).
$$
Then, $I_m\to\infty$ on a set of integers $m$ of natural density $1$.
\end{thm}

We follow the proof of Theorem \ref{thm1}.  Keeping
the notation introduced there, we see that
for large $J$,
$$
\PP \( \bigcap_{J\le j\le 3J/2} \overline{A_jB_jC_j} \) \le
\sum_{J\le j\le 3J/2} \frac{D^2}{j^2} + \prod_{J\le j\le 3
     J/2} \(1 - \frac{3}{4\sqrt{j}}\) \ll \frac{1}{J}\cdot
$$
For large $G$, define $J$ by $N_{J+1} < G \le N_{J+2}$.  Then
$G^{5/3} > N_{\fl{3J/2}+2}$ and $J\gg_M (\log G)/\log_2 G$.
Thus, for large~$G$,
$$
\PP\( \min_{G\leqslant N\leqslant G^{5/3}} \;
\max_{h(N)< n\le h(Nf_{N})} \frac{|S_n|}{s_n} \le K \) \ge 1 -
O\pfrac{1}{J} \ge 1 - O\pfrac{\log_2 G}{\log G}.
$$
The direct number theoretic analog of $|S_n|/s_n$ is
$$
\tw(m,t) := \frac{\big|\omega(m,t)-\sum_{p\le t} 1/p\big|}{\sqrt{\sum_{p\le t}
    (1-1/p)/p}}.
$$
By \eqref{Kubilius}, if $G$ is large and $G\le \sqrt{\log_2 x}$ (so
that $G^{5/3}f_{G^{5/3}} \le (\log_2 x)^{7/8}$),
then
$$
\frac{1}{x}\Big| \Big\{ m\le x : \min_{G\leqslant N\leqslant G^{5/3}} \;
\max_{h(N)<n\le h(Nf_N)} \tw(m,t) \le K \Big\} \Big| \ge 1 -
O\pfrac{\log_2 G}{\log G}.
$$
Since $\tw(m,t)= \varrho(m,t)+O\(1/\sqrt{\log_2 t}\,\)$, the first part of
the theorem follows.
\par\goodbreak
The second part is similar.  Note that
$\omega(n,x)-\omega\big(n,x^{1/\sqrt{\log_2 x}}\big) \le \sqrt{\log_2 x}$ for
$n\le x$, and, for brevity, write $g=g(\sqrt{x})$.
By \eqref{Kubilius} with $u:=\sqrt{\log_2 x}$, we have, for any fixed $K$
and large $x$,
\begin{align*}
\frac{1}{x} \biggl| \bigg\{ m\le x &: \min_{ \substack{N\ge g\\ Nf_N
\le \log_2 m }}
     \max_{N< \log_2 t \le Nf_N} \tw(m,t) \le K \bigg\} \biggr| \\
&\le \frac{1}{x} \left| \bigg\{
     \sqrt{x}\le m\le x : \min_{ \substack{N\ge g\\ Nf_N \le \LL(x) }}
     \max_{N\le \log_2 t\le Nf_N} \tw(t) \le K+2\bigg\}\right| +
\frac1{\sqrt{x}} \\
&\le \PP \( \inf_{\substack{N\ge g \\ Nf_N \le \LL(x)}}
      \max_{h(N)< n \le h(Nf_N)} \frac{|S_n|}{s_n} \le K+2 \) +
O\Big(\frac1{\log_2 x}\Big),
\end{align*}
where $\LL(x):=\log_2 x - \frac12 \log_3 x$.
Since $f_N\le N$, we have $N_{j+1}^* \le (N_j^*)^2$
in the notation of the proof
of Theorem~\ref{thm1}.  The interval $$\Big[(\log_2 x)^{1/10},
\LL(x)^{1/2}\Big]$$ therefore contains at least one interval $J_j$.  By
the proof of
Theorem \ref{thm1}, for large $x$,
the probability above does not exceed $\sum_{j\ge j_0} 1/j^2 \le
1/(j_0-1)$, where $j_0\to \infty$ as $x\to \infty$.

\noindent
{\bf Remarks.}  The upper bound $g^2$ of $N$ in the first part can
be sharpened.
By the same methods, similar results can be proved for a wide class of
additive arithmetic functions $r(m,t)=\sum_{p^a\| m} r(p^a)$ in
place of $\omega(m,t)$.

\noindent{\bf Acknowledgment.}  The authors are indebted to Walter
Philipp for helpful discussions on the use of almost sure
invariance principles.



\end{document}